# Parametric Design and Adaptive Sizing of Lattice Structures for 3D Additive Manufacturing

Jorge Manuel Mercado – Colmenero[1, 2] [0000-0001-8823-4896] Daniel Díaz – Perete [1] [0000-0002-2874-703X] Miguel Ángel Rubio– Paramio[1, 2] [0000-0003-1079-2323] and Cristina Martín – Doñate[1, 2] [0000-0002-3306-9456]

[1] Departament of Engineering Graphics Design and Projects, Campus Las Lagunillas s/n. A3 Building, Jaén 23071, Spain

[2] INGDISIG Jaén Research Group, Campus Las Lagunillas s/n. A3 Building, Jaén 23071, Spain

jmercado@ujaen.es

**Abstract.** The present research is developed into the realm of industrial design engineering and additive manufacturing by introducing a parametric design model and adaptive mechanical analysis for a new lattice structure, with a focus on 3D additive manufacturing of complex parts. Focusing on the landscape of complex parts additive manufacturing, this research proposes geometric parameterization, mechanical adaptive sizing, and numerical validation of a novel lattice structure to optimize the final printed part volume and mass, as well as its structural rigidity. The topology of the lattice structures exhibited pyramidal geometry. Complete parameterization of the lattice structure ensures that the known geometric parameters adjust to defined restrictions, enabling dynamic adaptability based on its load states and boundary conditions, thereby enhancing its mechanical performance. The core methodology integrates analytical automation with mechanical analysis by employing a model based in two-dimensional beam elements. The dimensioning of the lattice structure is analyzed using rigidity models of its sub-elements, providing an evaluation of its global structural behavior after applying the superposition principle. Numerical validation was performed to validate the proposed analytical model. This step ensures that the analytical model defined for dimensioning the lattice structure adjusts to its real mechanical behavior and allows its validation. The present manuscript aims to advance additive manufacturing methodologies by offering a systematic and adaptive approach to lattice structure design. Parametric and adaptive techniques foster new industrial design engineering methods, enabling the dynamic tailoring of lattice structures to meet their mechanical demands and enhance their overall efficiency and performance.

# 1 Introduction

Additive manufacturing technologies have witnessed a remarkable ascent in recent years, revolutionizing traditional manufacturing processes across various engineering areas. This transformative trajectory has been propelled by advancements in materials science, computational design, and manufacturing techniques, positioning additive manufacturing as a cornerstone technology with a wide versatility and applicability on each branch of engineering [1].

Among the different additive manufacturing techniques, MEX (Material Extrusion) stands out for its widespread adoption and versatility in manufacturing complex geometries and intricate structures. However, it is important to recognize that MEX presents certain challenges, particularly when fabricating beam lattice structures with small cross-sections, due to retraction issues and thermal-induced beam flexion during layer deposition. While this study focuses on MEX, future research should also explore other additive manufacturing techniques, such as material jetting, binder jetting, or powder bed fusion, which may offer better precision for such complex geometries [2]. MEX technology, commonly known as FDM (Fused Deposition Modeling) or FFF (Fused Filament Fabrication), offers a cost-effective and accessible means to translate digital designs into physical objects with high precision and fidelity [3-7].

MEX technology offers flexibility in fabricating complex geometries and free-form surfaces, although its capability to fabricate lattice-type structures depends significantly on the specific lattice geometry. While 2.5D lattice structures are relatively easy to manufacture and are often incorporated into slicing software, more complex 3D beam, TPMS, and stochastic lattice structures present challenges for MEX due to issues like retraction, beam deflection, and resolution limitations [8, 9]. In the present manuscript, we focus on exploring the manufacturability of a pyramidal lattice structure using MEX 3D additive manufacturing technology, recognizing its potential limitations. Although references [10-14] involve lattice structures produced by other AM (Additive Manufacturing) technologies, they provide a useful comparison for understanding the broader context of lattice structure fabrication and mechanical properties.

In the realm of additive manufacturing, several lattice geometries have emerged as seminal designs, each with distinct characteristics and applications. Among these, the diamond, gyroid, and octet lattice topologies have garnered considerable attention for their exceptional mechanical properties and geometric efficiency [15-17]. These lattice geometries serve as inspirations for the development and optimization of novel lattice structures adapted to specific engineering requirements.

In this context, the present manuscript proposes the design of a novel lattice structure characterized by pyramidal geometry and adaptive parameterization. By harnessing the inherent properties of pyramidal elements, our proposed lattice design offers unparalleled versatility in adapting to diverse filling requirements within industrial components, without the need for auxiliary elements that facilitate its manufacturing, such as internal supports. Furthermore, the defined lattice structure design is completely parameterized which facilitates precise dimensioning and customization to adapt varying mechanical and geometric constraints, thereby enhancing overall performance and manufacturability.

In summary, the methodology described in the present manuscript aims to explore the design, optimization, and application of lattice structures within the realm of 3D additive manufacturing.

Through a comprehensive review of existing lattice geometries and the introduction of our novel pyramidal lattice design, it is seeked to advance the capabilities and potential of additive manufacturing technologies in industrial design engineering applications.

## 2 Lattice parameterized geometric design

This section describes the geometric characteristics and parameterization process of the lattice-type structure proposed in the present manuscript. Figure 1 shows the set of designed lattice-type structures as well as the geometric definition of their lattice unit cells. As can be seen, the lattice unit cell is designed as a pyramidal geometry. The layout of each edge that composes it allows its manufacturing using 3D additive manufacturing technology to be performed without the need to include any auxiliary support type element. This lattice unit cell design confers greater mechanical rigidity to the structural assembly, owing to its pyramidal topology. This geometric shape is known to have high structural rigidity. This is mainly because of its ability to uniformly distribute the set of forces applied to it and its ability to resist deformations. Likewise, its shape is composed of four triangular faces that connect to each other at shared points, creating a solid three-dimensional structure. Furthermore, this pyramidal geometry efficiently uses the materials used. It minimizes its volume without promising its resistance, which makes it ideal for applications where a high strength-to-weight ratio is required.

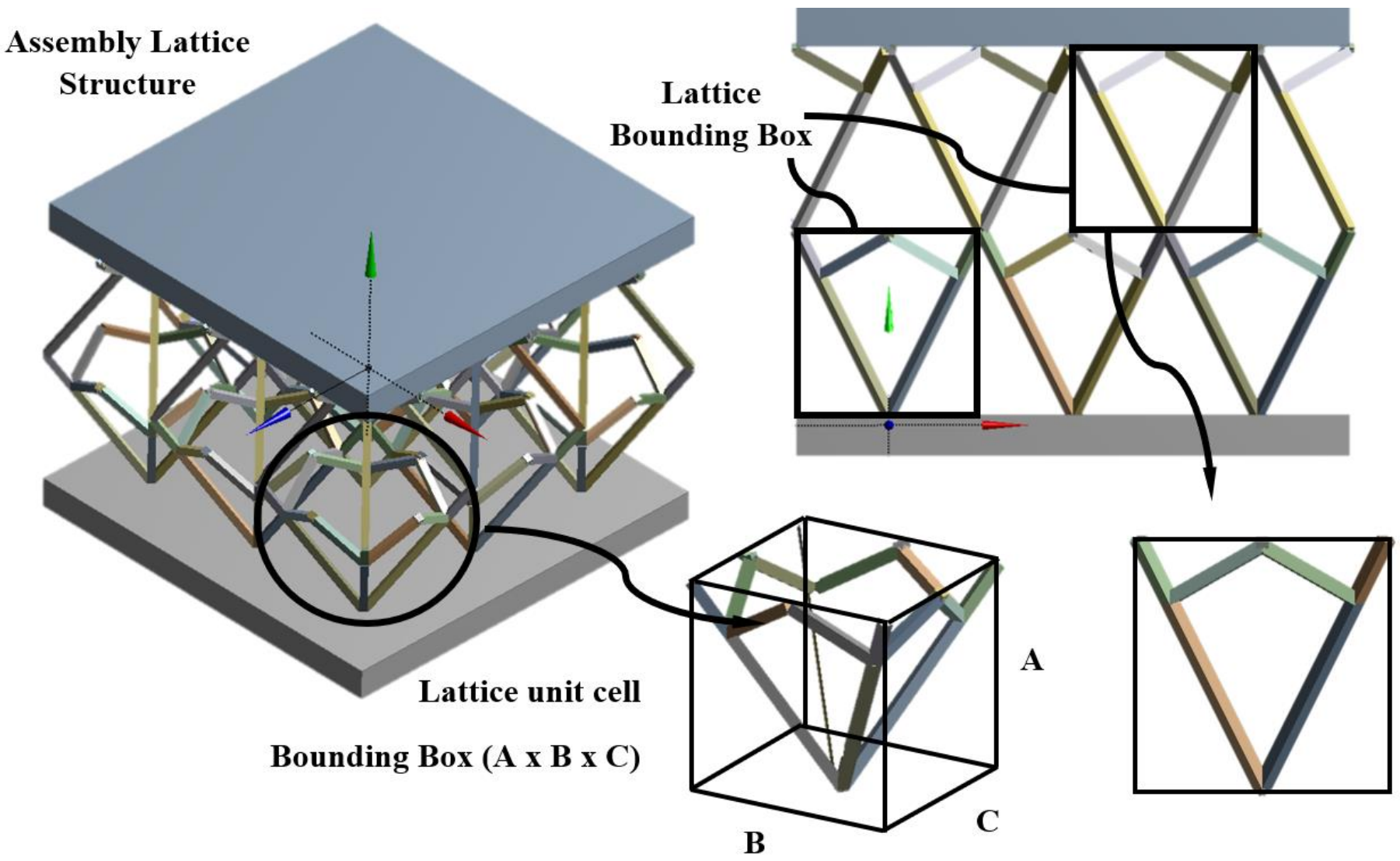

**Fig. 1.** Geometric design of the proposed lattice structure.

The parameterization process of the lattice unit cell was developed based on the dimensions of its thickness and the size of the edge that establishes its Bounding Box. As shown in Equation 1 and Figure 2, from these geometric parameters and together with the definition of aspect ratios, the dimensions of the final Bounding Box that contains the geometry of the lattice unit cell are completed.

$$\begin{Bmatrix} A \\ B \\ C \end{Bmatrix} = \begin{Bmatrix} \gamma_x \\ \gamma_y \\ \gamma_z \end{Bmatrix} \cdot L \qquad (1)$$

Where L [mm] represents the dimension of the edge of the Bounding Box of the lattice unit cell and $\gamma_x$, $\gamma_y$, $\gamma_z$ represent the aspect ratios that affect the dimensions of the Bounding Box of the lattice unit cell. Next, establishing as independent variables the definitive dimensions of the Bounding Box of the lattice unit cell, A (height), B (width) and C (depth) [mm] (see Figure 1 and Figure 2), and the cross section thickness of the lattice unit cell bars, T [mm], the rest of the geometric parameters that determine the topology and structure of the lattice unit cell are established, see Figure 2, Equation 2, Equation 3 and Equation 4.

$$D = \sqrt{\left(B/2\right)^2 + \left(C/2\right)^2 + A^2} \qquad (2)$$

$$\beta/2 = \mathrm{ArcSin}\left(\frac{B/2}{D}\right) \qquad (3)$$

$$E = B/2 \cdot \mathrm{Sin}\left(\pi/2 - \beta/2\right) \; ; \; F = B/2 \cdot \mathrm{Cos}\left(\pi/2 - \beta/2\right) \qquad (4)$$

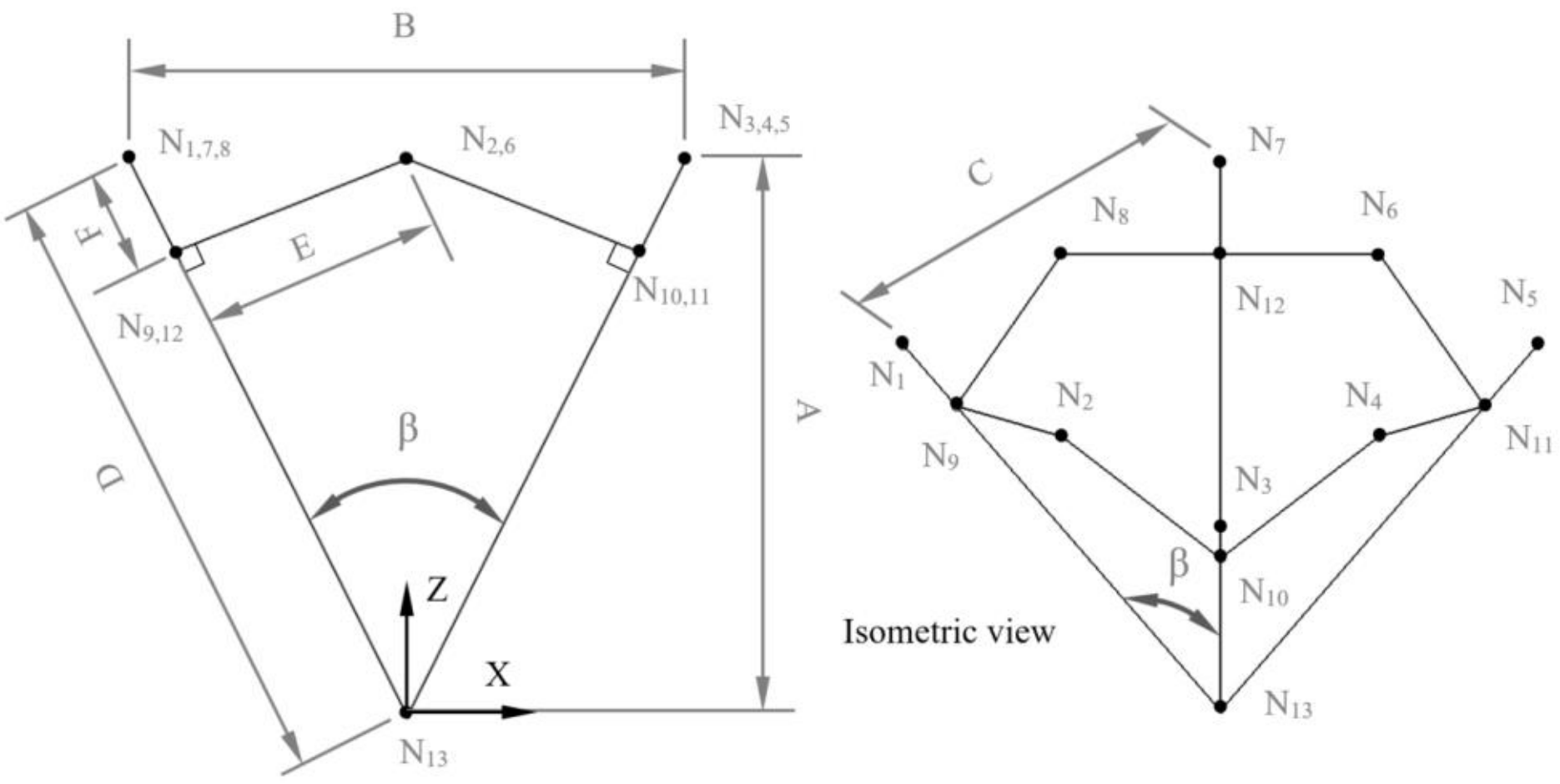

**Fig. 2.** Parametric design of the proposed lattice unit cell.

The manufacturability of the proposed lattice unit cell design is validated using additive manufacturing software. These analyses confirm that the lattice design does not require support structures during printing, and it can be manufactured without difficulty using 3D additive technology. The proposed design is optimized to eliminate overhangs and complex geometries that typically necessitate support materials, ensuring a straightforward and efficient printing process.

Figure 2 shows the geometric distribution of each unidirectional element that makes up the proposed lattice unit cell. Moreover, it also shows the geometric distribution of each of the nodes that make up the proposed lattice geometry. To complete the geometric parameterization process of this design, Table 1 shows the expression of the Cartesian coordinates of each node from the geometric parameterization expressions presented in Equation 2, Equation 3 and Equation 4.

**Table 1.** Parametric definition of the nodes coordinates of the proposed lattice unit cell.

| Lattice Nodes | Cartesian coordinates [mm] |
|---|---|
| 1, 3 | [± B/2, – C/2, A] |
| 5, 7 | [± B/2, C/2, A] |
| 2, 6 | [0, ± C/2, A] |
| 4, 8 | [± B/2, 0, A] |
| 9, 12 | [–B/2 + F·Sin (β/2), ± C/2 ± F·Sin (β/2), A– F·Cos (β/2)] |
| 10, 11 | [B/2 – F·Sin (β/2), ± C/2 ± F·Sin (β/2), A– F·Cos (β/2)] |
| 13 | [0 , 0, 0] |

After describing the geometry of the lattice unit cell and its parameterization process, we proceeded to analytically evaluate its structural behavior. To do this, as shown in Figure 3, the rigidity model is applied, considering that the elements that make up the lattice unit cell are the unidirectional 3D beam type, with a constant cross section and three degrees of freedom in each of its extreme nodes.

The 3D beam element is selected for this numerical analysis due to its ability to capture combined blending and axial effects in slender structures, making it an appropriate choice for lattice members with high aspect ratios. This element type also provides a good balance between computational efficiency and accuracy in predicting the global and local deformation behavior of the lattice unit cell.

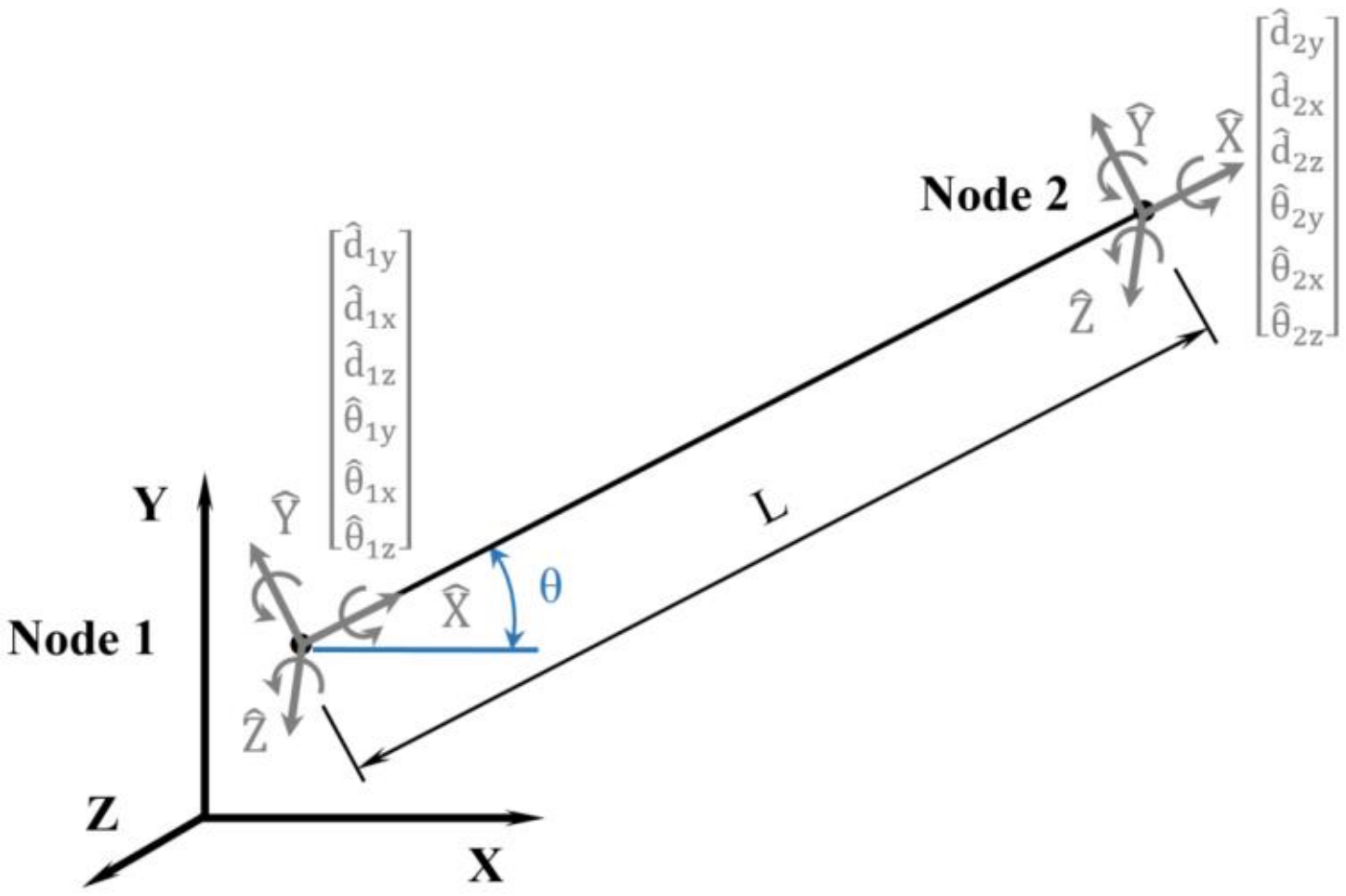

**Fig. 3.** Definition of the 3D unidirectional beam type element.

First, as established in Equation 5, the transformation matrix $\underline{T}$ of the defined beam element, which relates the local displacement field $\underline{\hat{d}}$ to the global reference system and the global displacement vector $\underline{d}$, is determined.

$$\underline{\hat{d}} = \underline{T} \cdot \underline{d} = \begin{Bmatrix} \hat{d}_{1x} \\ \hat{d}_{1y} \\ \hat{d}_{1z} \\ \hat{\theta}_{1x} \\ \hat{\theta}_{1y} \\ \hat{\theta}_{1z} \\ \hat{d}_{2x} \\ \hat{d}_{2y} \\ \hat{d}_{2z} \\ \hat{\theta}_{2x} \\ \hat{\theta}_{2y} \\ \hat{\theta}_{2z} \end{Bmatrix} = \begin{bmatrix} C_x & C_y & C_z & 0 & 0 & 0 & 0 & 0 & 0 \\ S_x & S_y & S_z & 0 & 0 & 0 & 0 & 0 & 0 \\ 0 & 0 & 0 & 1 & 0 & 0 & 0 & 0 & 0 \\ 0 & 0 & 0 & C_x & C_y & C_z & 0 & 0 & 0 \\ 0 & 0 & 0 & S_x & S_y & S_z & 0 & 0 & 0 \\ 0 & 0 & 0 & 0 & 0 & 0 & C_x & C_y & C_z \\ 0 & 0 & 0 & 0 & 0 & 0 & S_x & S_y & S_z \end{bmatrix}_{sym} \cdot \begin{Bmatrix} d_{1x} \\ d_{1y} \\ d_{1z} \\ \theta_{1x} \\ \theta_{1y} \\ \theta_{1z} \\ d_{2x} \\ d_{2y} \\ d_{2z} \\ \theta_{2x} \\ \theta_{2y} \\ \theta_{2z} \end{Bmatrix} ; \begin{matrix} C_x = \mathrm{Cos}\,(\theta_x) \\ C_y = \mathrm{Cos}\,(\theta_y) \\ C_z = \mathrm{Cos}\,(\theta_z) \\ S_x = \mathrm{Sin}\,(\theta_x) \\ S_y = \mathrm{Sin}\,(\theta_y) \\ S_z = \mathrm{Sin}\,(\theta_z) \end{matrix} \quad (5)$$

After determining the transformation matrix of the 3D beam element, its stiffness matrix $\underline{K}$ is defined based on static equilibrium equations and the 3D beam geometric and material properties, as shown in Equation 6.

$$\underline{K} = \frac{E}{L} \cdot \begin{bmatrix} A/L & 0 & 0 & 0 & 0 & 0 & -A/L & 0 & 0 & 0 & 0 & 0 \\ 0 & 12 \cdot I_z/L^3 & 0 & 0 & 0 & 6 \cdot I_z/L^2 & 0 & -12 \cdot I_z/L^3 & 0 & 0 & 0 & 6 \cdot I_z/L^2 \\ 0 & 0 & 12 \cdot I_y/L^3 & 0 & -6 \cdot I_y/L^2 & 0 & 0 & 0 & -12 \cdot I_y/L^3 & 0 & -6 \cdot I_y/L^2 & 0 \\ 0 & 0 & 0 & G \cdot J/L & 0 & 0 & 0 & 0 & 0 & -GJ/L & 0 & 0 \\ 0 & 0 & -6 \cdot I_y/L^2 & 0 & 4 \cdot I_y/L & 0 & 0 & 0 & 6 \cdot I_y/L^2 & 0 & 2 \cdot I_y/L & 0 \\ 0 & 6 \cdot I_z/L^2 & 0 & 0 & 0 & 4 \cdot I_z/L & 0 & -6 \cdot I_z/L^2 & 0 & 0 & 0 & 2 \cdot I_z/L \\ -A/L & 0 & 0 & 0 & 0 & 0 & A/L & 0 & 0 & 0 & 0 & 0 \\ 0 & -12 \cdot I_z/L^3 & 0 & 0 & 0 & -6 \cdot I_z/L^2 & 0 & 12 \cdot I_z/L^3 & 0 & 0 & 0 & -6I_z/L^2 \\ 0 & 0 & -12 \cdot I_y/L^3 & 0 & 6 \cdot I_y/L^2 & 0 & 0 & 0 & 12 \cdot I_y/L^3 & 0 & 6 \cdot I_y/L^2 & 0 \\ 0 & 0 & 0 & 0 & 0 & 0 & 0 & 0 & 0 & GJ/L & 0 & 0 \\ 0 & 0 & -6 \cdot I_y/L^2 & G \cdot J/L & 2 \cdot I_y/L & 0 & 0 & 0 & 6 \cdot I_y/L^2 & 0 & 4 \cdot I_y/L & 0 \\ 0 & 6 \cdot I_z/L^2 & 0 & 0 & 0 & 2 \cdot I_z/L & 0 & -6 \cdot I_z/L^2 & 0 & 0 & 0 & 4 \cdot I_z/L \end{bmatrix}_{SYM} \quad (6)$$

Where I [$mm^4$] represents the inertia of the cross section of the beam element, E [MPa] represents the elastic modulus of the material associated with the beam element, A [$mm^2$] represents the constant cross section area of the beam-type element analyzed, and L [mm] is the length of the beam element, as shown in Figure 3. Next, from Equations 5 and 6, a system of equilibrium equations of a beam-type element is established for the global reference system. Equation 7 determines both the global displacement field of the 3D beam-type element as well as the reaction forces and moments at its extreme nodes.

$$\underline{f} = \underline{T}^{T} \cdot \underline{K} \cdot \underline{T} \cdot \underline{d} \tag{7}$$

Where $\underline{f}$ represents the global force vector, including axial forces, shear forces, and bending/torsional moments, $\underline{K}$ represents the stiffness matrix in the global coordinate system and $\underline{d}$ represents the displacement field (translations and rotations in the global reference system). Finally, to validate the design process of the lattice unit cell proposed in this study, its structural behavior is evaluated. On the one hand, the rigidity model and the static equilibrium equations, see Equation 7, are applied to obtain an analytical solution of the structural behavior of the lattice unit cell. While the analytical model provides a reliable prediction of the structural behavior of the lattice unit cell, a numerical analysis using FEM was performed to further validate the results. FEM validation is widely used in engineering to ensure that analytical solutions hold under a variety of practical boundary conditions and load scenarios. By confirming the agreement between the analytical and numerical approaches, this validation process ensures the robustness and applicability of the proposed lattice design in real-world applications, particularly when facing more complex conditions not easily addressed by purely analytical models. The numerical analysis was performed using the commercial FEM software ANSYS Mechanical (ANSYS, Inc., Canonsburg, PA, USA) [18]. Furthermore, during the definition of the preprocessing phase of this numerical simulation, the element type defined to model the set of the lattice unit cell is a second-order unidirectional beam type and the element size established was 5 mm, see Figure 4.

Table 2 shows the geometric dimensions, mechanical and elastic properties of the plastic material, type of element and mesh size, boundary conditions, and defined load scenario (see Figure 4) used to evaluate the mechanical and structural behavior of the lattice unit cell.

**Table 2.** Lattice unit cell structural analysis setup.

| | Parameter | Unit | Value |
|---|---|---|---|
| Lattice geometry | Bounding Box dimensions [A, B, C] | mm | 100 |
| | Thickness | mm | 4.0 |
| | $\gamma_x$, $\gamma_y$, $\gamma_z$ | – | 1 |
| | Cross section geometry | – | Square |
| Plastic Material | Trade name | – | PETG |
| | Elastic Modulus, E | MPa | 2,800 |
| | Poisson’s ratio, υ | – | 0.33 |
| Mesh | Element size | mm | 5.0 |
| | Element type | – | Beam |
| Boundary condition | Support | – | Fixed |
| Load scenario | Force | N | 100 |

A square cross-section is selected for the lattice struts to facilitate the manufacturability of the lattice structure using 3D printing technology. This cross-section ensures better control of layer deposition during the printing process, providing consistent wall thickness and enhanced mechanical stability. Although circular cross-sections are more common in lattice structures, the square cross-section was chosen here to simplify the manufacturing process without compromising structural performance.

Figure 4 shows the boundary conditions and load scenario to which the lattice unit cell is subjected, as well as the resulting displacement field. As can be seen, the maximum displacement obtained, located at the application node of the load scenario, is equal to 0.0714 mm. Moreover, this numerical result is analogous to the analytical result obtained after applying the rigidity model and static equilibrium equations (see Equation 7) to the lattice unit cell. Thus, the analytical model used is validated, as well as the parameterization process proposed for the design of the present lattice.

This validation confirms that the Cartesian coordinates used in the definition of the nodes for the numerical and analytical mechanical and structural analysis of the unit cells proposed in the present manuscript are valid.

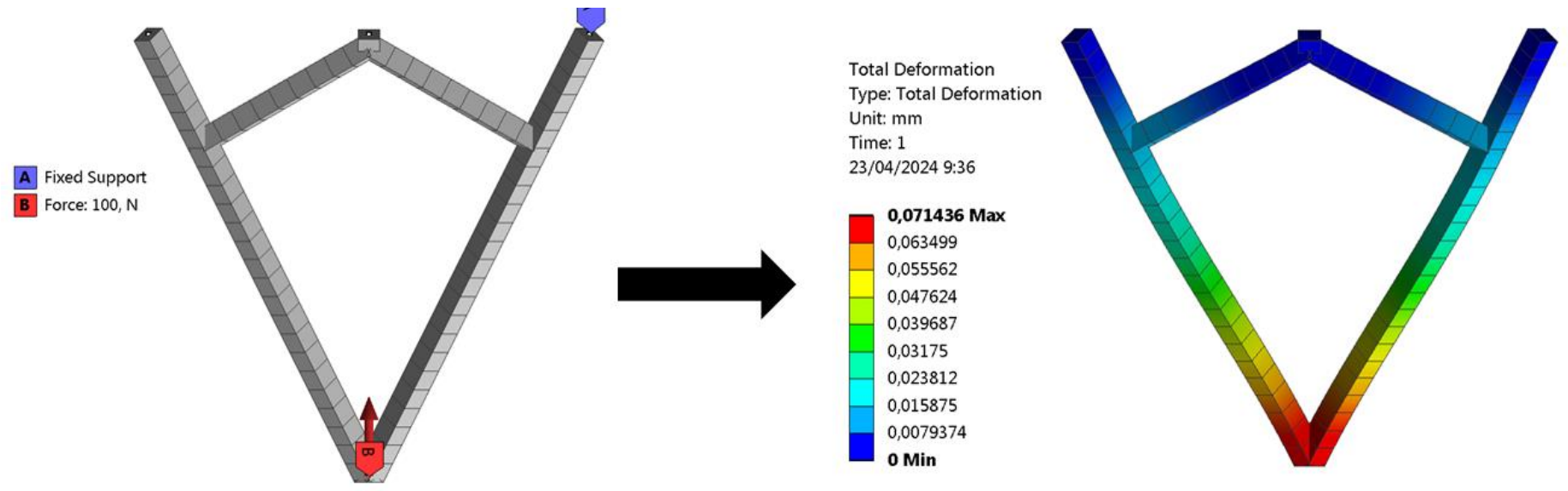

**Fig. 4.** Boundary conditions and load scenario for the lattice unit cell FEM analysis.

Finally, the process of parameterizing the geometry of the lattice unit cell, together with the analytical model that determines its mechanical behavior and structure, allows the determination of its associated structural rigidity. Therefore, by applying the principle of superposition, the structural behavior of a matrix or set of lattices proposed in the present manuscript can be evaluated analytically.

# 3 Results and discussion

Next, to evaluate the structural behavior of the proposed lattice design, a hollow geometry was selected as a case study (see Figure 5 and Figure 6). The geometry is contained within a bounding box with dimensions of 93.6 mm x 167 mm x 13 mm. All surfaces of the hollow geometry have a uniform wall thickness of 1.5 mm. This setup represents a realistic industrial component where the internal lattice structures occupy the interior volume, optimizing weight while maintaining structural integrity. The dimensions of the bounding box and the consistent wall thickness ensure uniformity in the comparison of the different lattice structures.

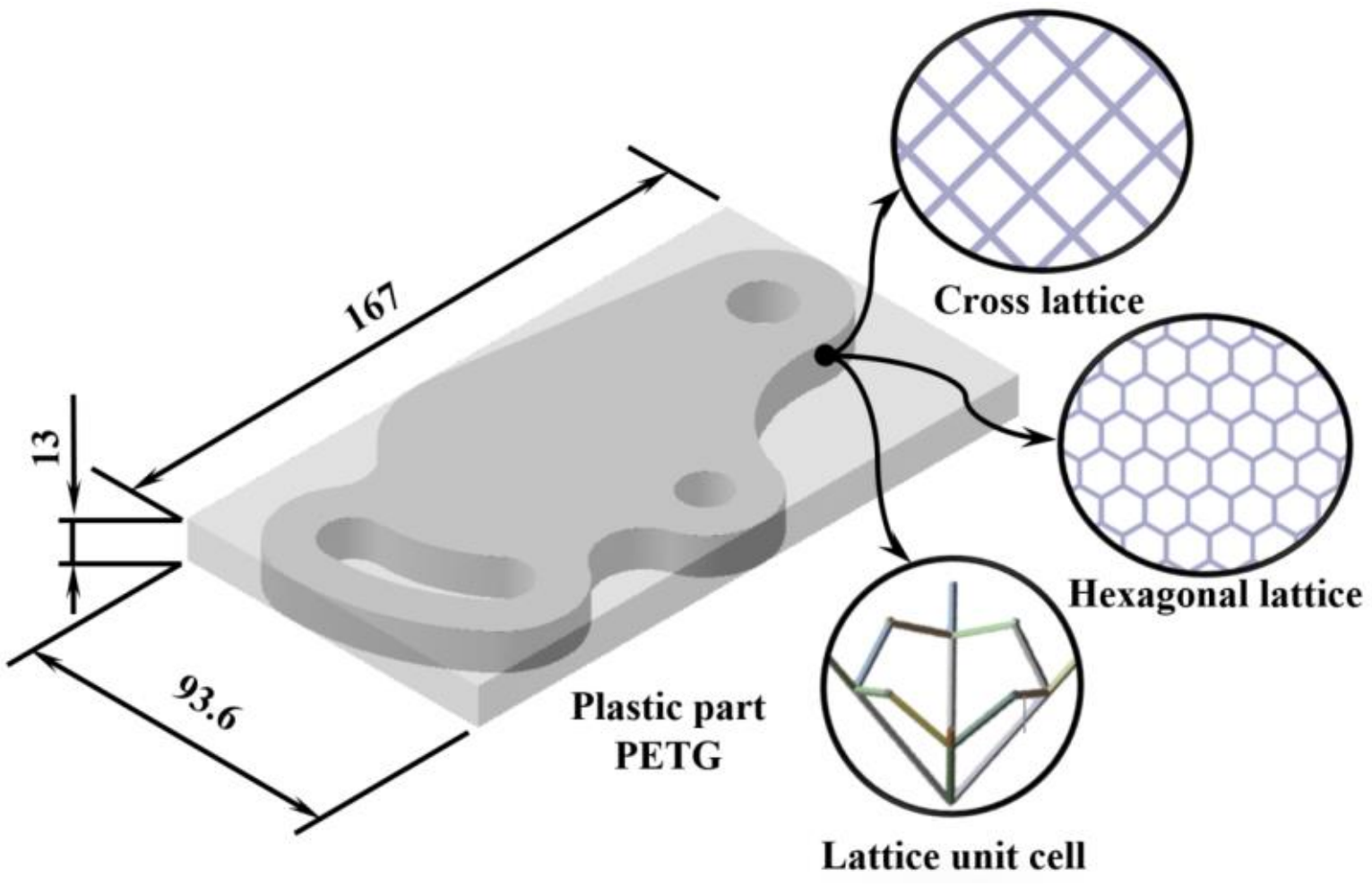

**Fig. 5.** Definition of the plastic part case study for the proposed lattice structure validation.

To complete the interior volume of the geometry under study, together with the lattice-type structure proposed in this manuscript, two lattice-type patterns commonly used in 3D additive manufacturing processes were defined, namely, the cross and hexagonal patterns. The selected comparison lattice structures, extruded along the Z – axis, are widely used in commercial additive manufacturing processes, particularly in MEX technology. These 2.5D lattice structures represent common benchmarks in industrial applications, allowing for a practical evaluation of the proposed lattice unit cell. The comparison aims to demonstrate how the new lattice design, with its 3D beam configuration, can outperform well-established structures in terms of material optimization and mechanical performance. This approach ensures the industrial relevance of the study, as these 2.5D lattices are widely adopted in various applications, making them suitable reference points for this analysis. For all lattice structures analyzed in the case study, including the proposed lattice unit cell, the hexagonal lattice, and the cross lattice, the thickness of the walls is consistently set to 0.5 mm. This uniform thickness ensures a fair comparison of the structural performance and optimization of the material used in additive manufacturing across the different lattice geometries.

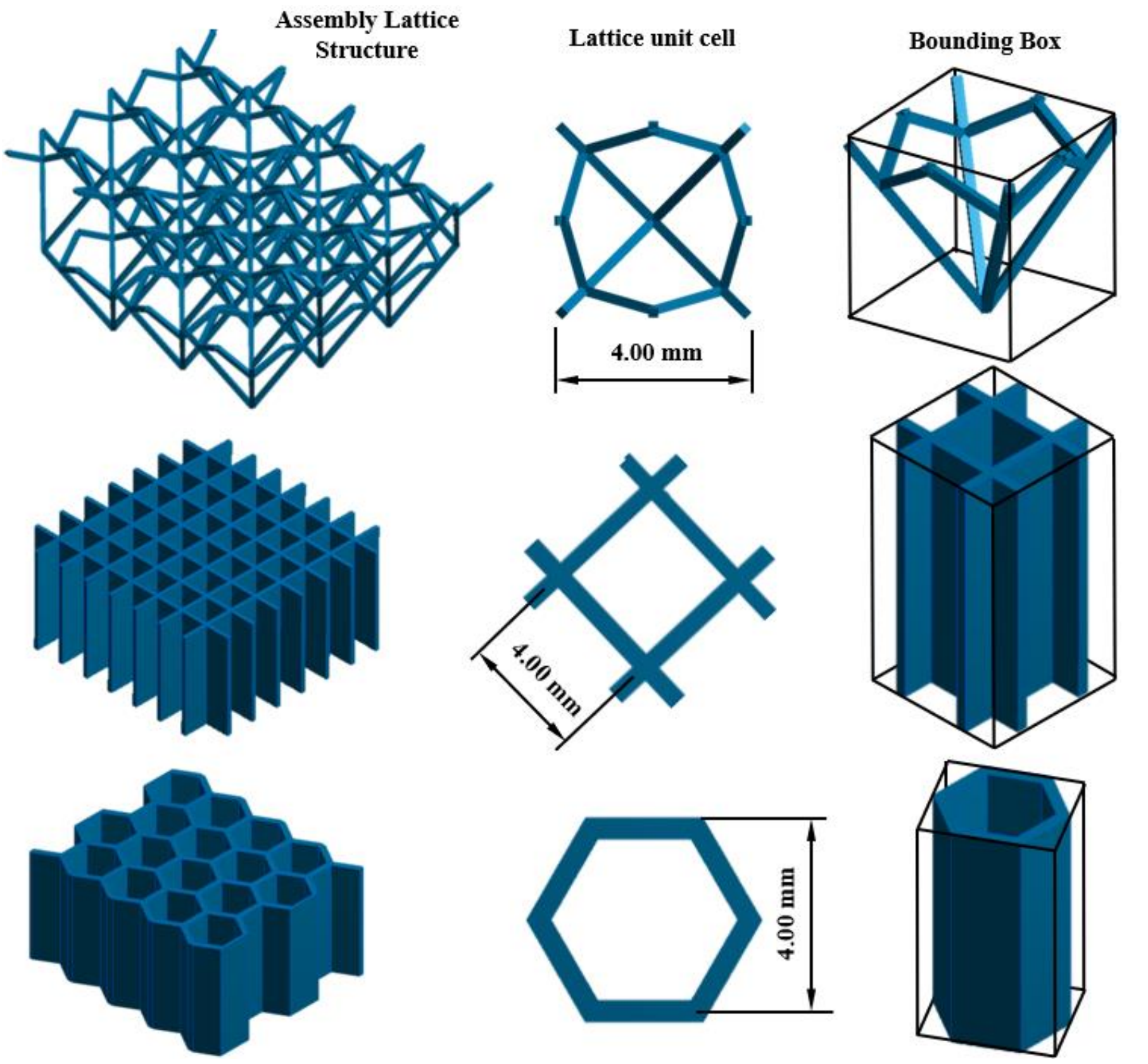

**Fig. 6.** Design and geometric characteristics of the reticular structures used for the validation

Thus, the structural behavior of the lattice unit cell proposed in this study can be verified and validated against reticular structures frequently used in 3D additive manufacturing processes. To do this, the commercial FEM software ANSYS Mechanical [18] is used to perform numerical analyses of the geometry plastic part case study, together with the reticular structures under study. Figure 6 shows the loading scenarios and boundary conditions defined for each numerical simulation. In particular, three different load scenarios were established, as the objective is to evaluate the structural behavior of the proposed case study in the principal axis X, Y and Z. The geometry meshing process (see Figure 7) is performed using second-order tetrahedral elements with a size of 0.5 mm for the lattice structures (cross, hexagonal and proposed lattice unit cell) and 3 mm for the surface of the plastic part case study. Moreover, the meshing process is performed automatically based on the type and size of the element used for the mesh, both for the inner region and the outer region. These element sizes are selected based on a mesh convergence analysis, which demonstrated that further refinement of the mesh produced negligible changes in displacement and stress values.

This ensures that the selected mesh provides sufficient accuracy while maintaining computational efficiency. Furthermore, the interior lattice elements and the exterior solid were structurally connected using node-to-node coupling at their interface, ensuring that the degrees of freedom are properly transferred between the interior and exterior elements, preserving structural continuity in the simulation. The plastic material used is PETG, and its elastic and mechanical properties are defined in Table 2.

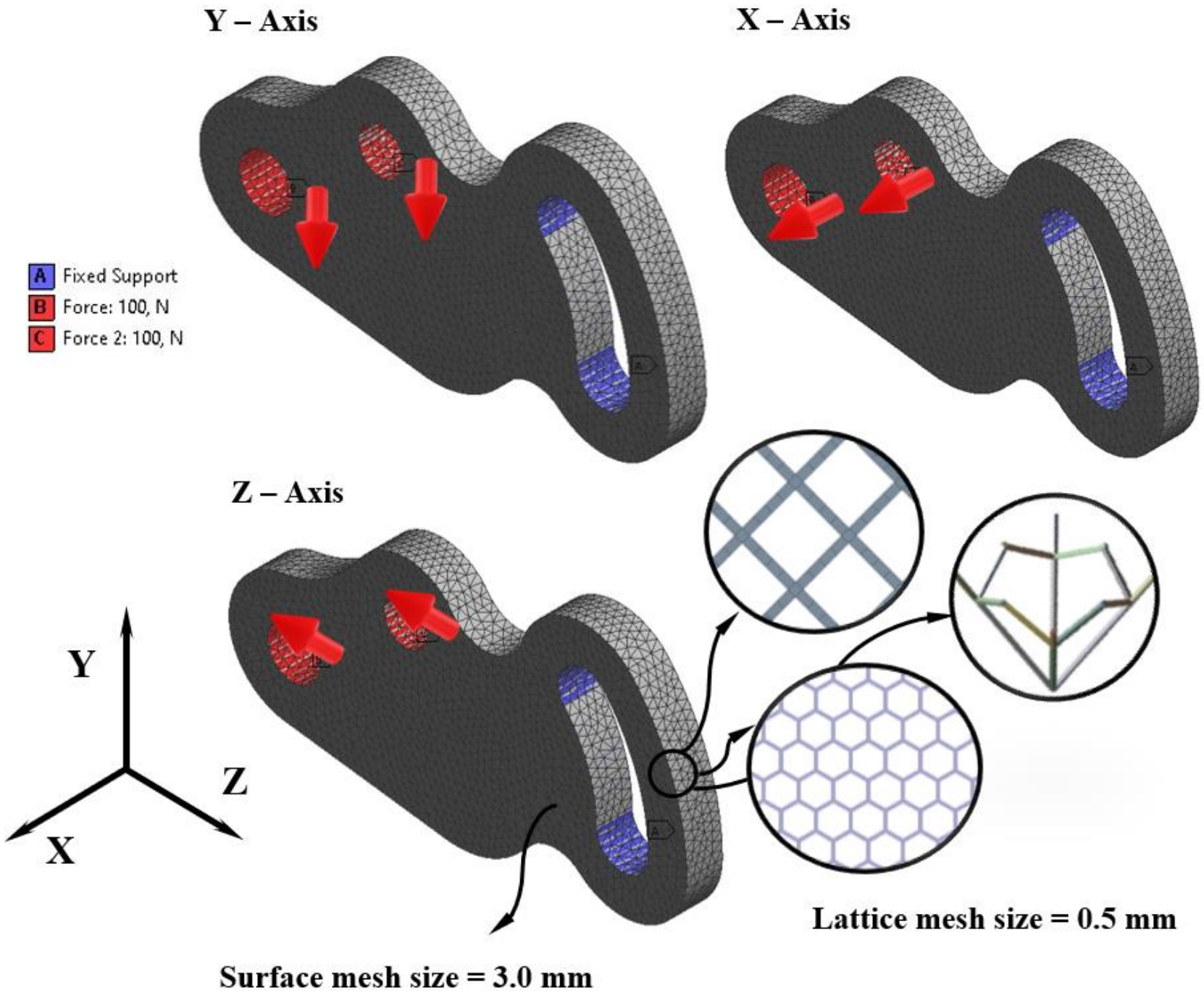

**Fig. 7.** Load scenarios and boundary conditions for the case study analysis.

Table 3, Figure 8, Figure 9 and Figure 10 show the maximum value of the total displacement field and von Mises stress field for each reticular structure analyzed and for each defined load scenario. As can be seen, the reticular structure that presents the lowest maximum displacement for each scenario analyzed is the cross-lattice type. However, the magnitude of this parameter does not determine the general structural performance of each reticular structure, or the volume and mass performance that are present in each structural scenario analyzed.

To do this, and from the results obtained and the mass of each geometry analyzed, the value of the structural efficiency parameter is determined (see Table 3). This parameter is determined as the ratio between the maximum displacements obtained and the mass of each geometry. Thus, using this parameter, the reticular structure that best optimizes the plastic material used in relation to its structural rigidity can be verified.

As can be seen in Table 3, the lattice proposed in the present manuscript is the one that obtains the best structural efficiency because it is the one that best optimizes the volume of plastic material necessary for its manufacture. Furthermore, the parameter improvement represents the relative difference between the magnitude of structural efficiency for each reticular geometry analyzed. It can be seen how the design of the proposed lattice unit cell improves the results obtained for the cross and hexagonal patterns, especially in those mechanical stresses or load scenarios that result in a greater magnitude in the maximum value of the displacement field. For example, in the Y and Z load application directions, the percentage of structural efficiency improvement for the proposed lattice design is between 32~39%, compared to the cross and hexagonal reticular structures.

This finding aligns with recent studies on lattice structures and their applications in additive manufacturing. For instance, Ye et al. (2021) and Naboni et al. (2020) highlight the importance of lattice topologies in optimizing strength and energy absorption, particularly in low-density structures, which is consistent with our results in terms of the proposed lattice's superior material efficiency [19, 20]. Additionally, Khiavi et al. (2022) show that lattice structures with optimized structure based topologies exhibit enhanced mechanical properties, supporting our finding that the proposed lattice design achieves a favorable balance between rigidity and material efficiency [21]. Furthermore, Graziosi et al. (2022) emphasize the role of beam-based lattice structures in energy-absorbing applications, which also resonates with the structural behavior observed in our lattice design under various loading conditions [22]. The comparison of our results with these studies confirms the mechanical robustness of our lattice structure and positions it as a competitive design for applications requiring high strength to weight ratios and structural efficiency. This also reflects the relevance of our lattice design for energy absorbing applications, as noted by the significant performance gains in load-bearing scenarios.

**Table 3.** Numerical results obtained for each lattice geometry analyzed.

| Load direction | Geometry | Max. Displacement [mm] | Max. Stress [MPa] | Mass [kg] | Structural efficiency [mm/kg] | Improvement [%] |
|---|---|---|---|---|---|---|
| Y – Axis | Proposed lattice | 0.518 | 18.904 | $68.493 \cdot 10^{-3}$ | 7.563 | – |
| | Cross lattice | 0.400 | 12.650 | $88.022 \cdot 10^{-3}$ | 4.544 | 39.918 |
| | Hexagonal lattice | 0.437 | 10.295 | $83.823 \cdot 10^{-3}$ | 4.953 | 34.510 |
| X – Axis | Proposed lattice | 0.045 | 5.644 | $68.493 \cdot 10^{-3}$ | 0.432 | – |
| | Cross lattice | 0.035 | 3.533 | $88.022 \cdot 10^{-3}$ | 0.398 | 9.028 |
| | Hexagonal lattice | 0.035 | 3.393 | $83.823 \cdot 10^{-3}$ | 0.418 | 3.349 |
| Z – Axis | Proposed lattice | 4.865 | 46.689 | $68.493 \cdot 10^{-3}$ | 71.029 | – |
| | Cross lattice | 3.956 | 112.070 | $88.022 \cdot 10^{-3}$ | 44.943 | 36.726 |
| | Hexagonal lattice | 4.011 | 59.323 | $83.823 \cdot 10^{-3}$ | 47.851 | 32.632 |

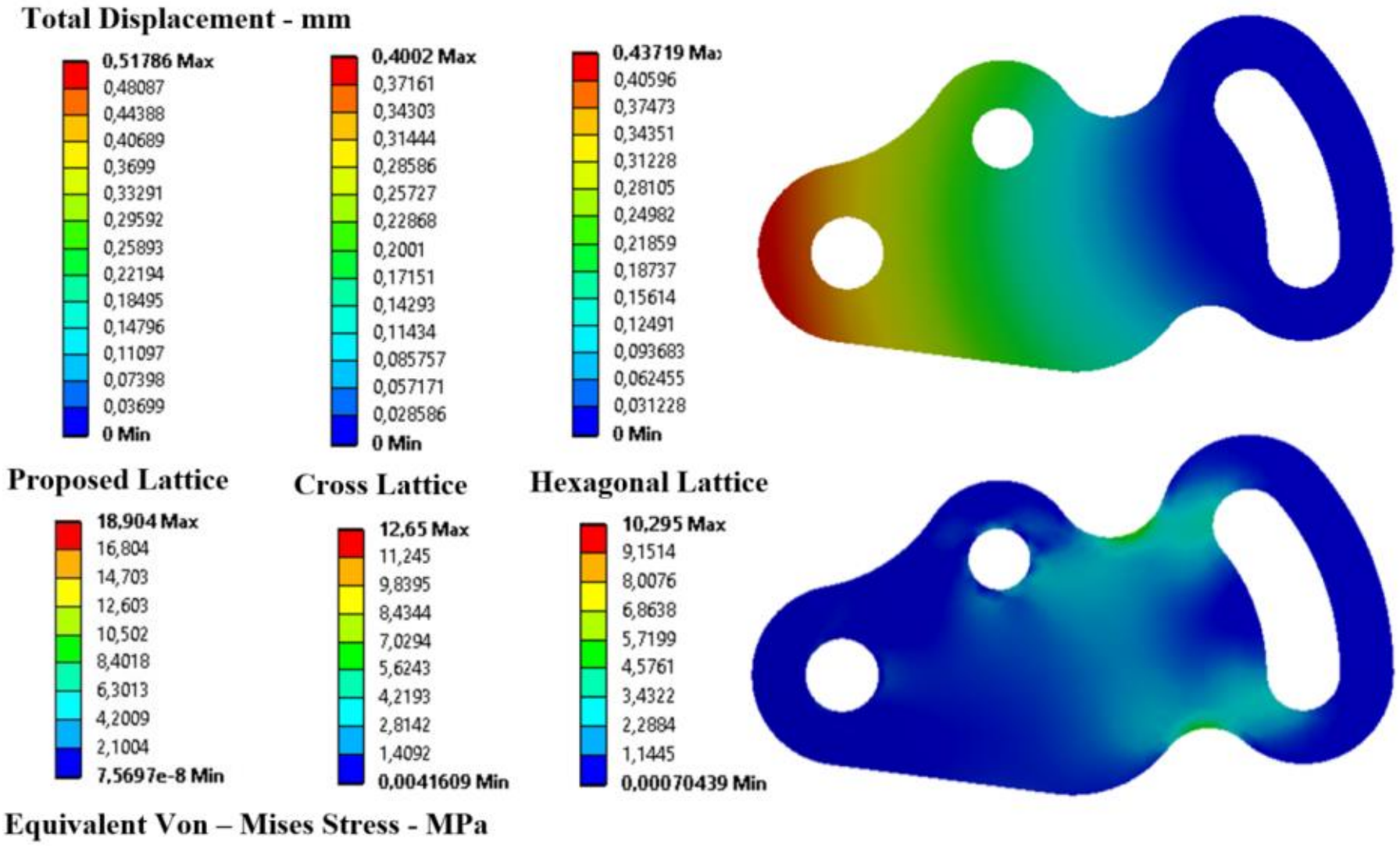

**Fig. 8.** Numerical results for each lattice geometry analyzed, Y–Axis mechanical study.

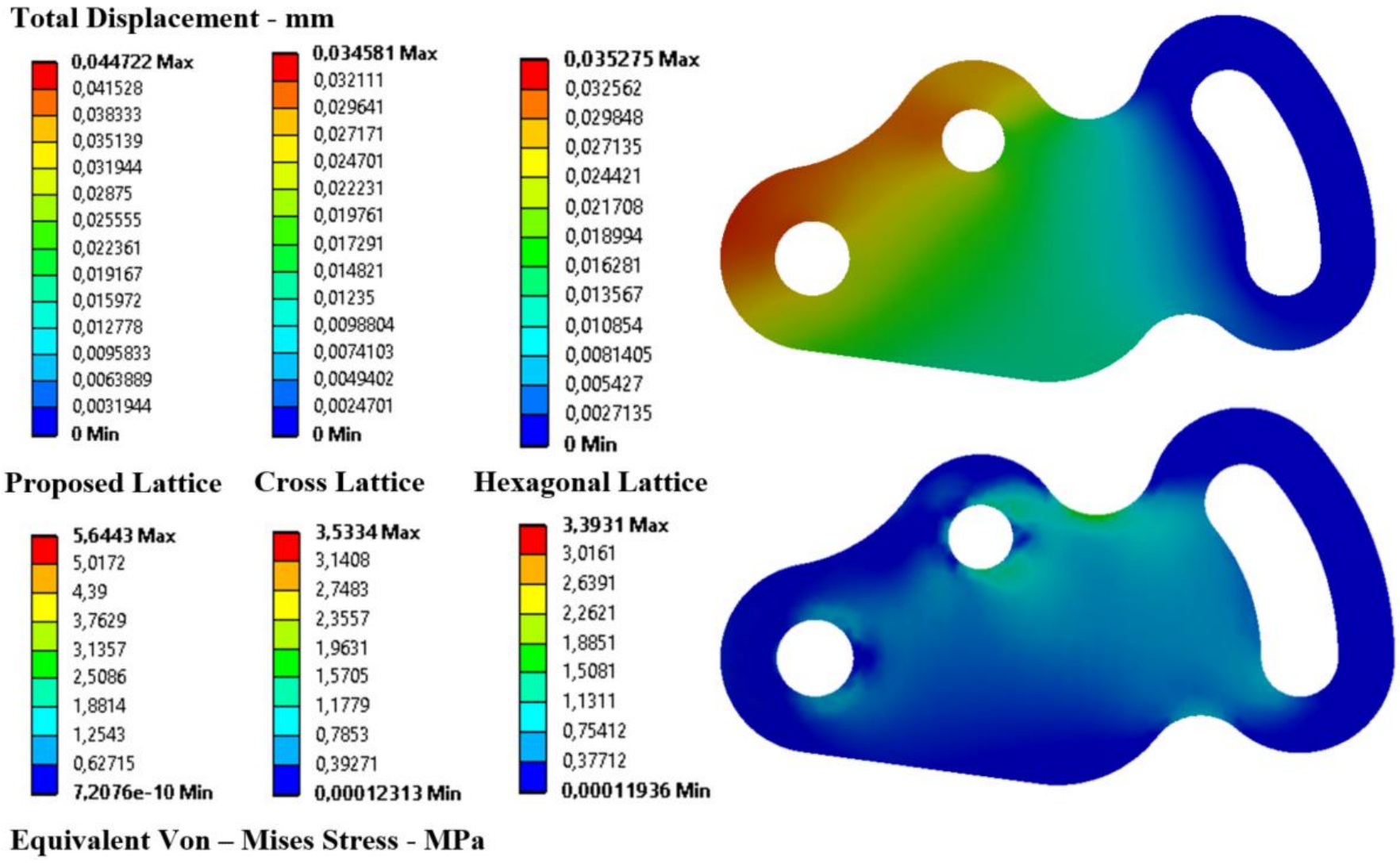

**Fig. 9.** Numerical results for each lattice geometry analyzed, X–Axis mechanical study.

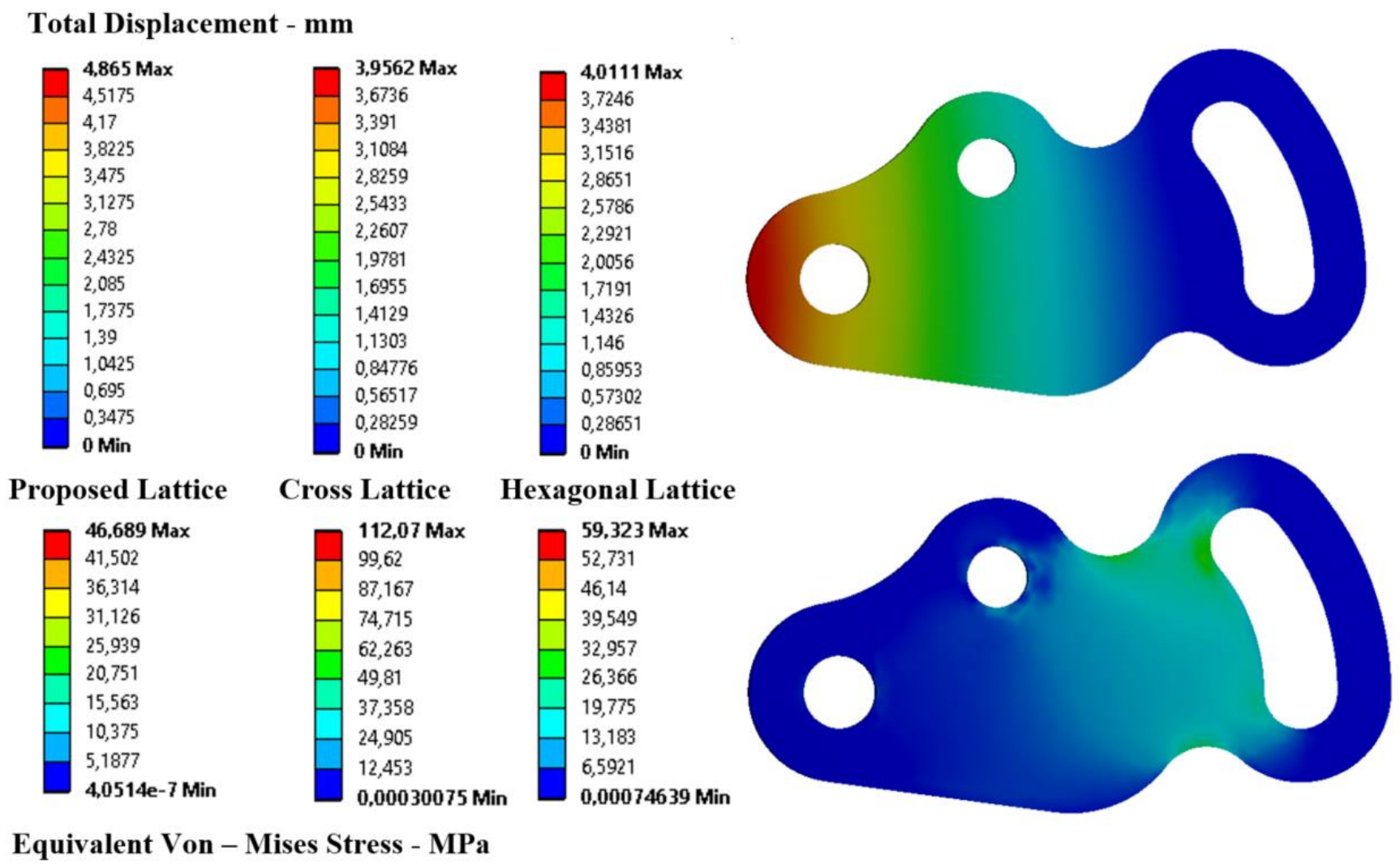

**Fig. 10.** Numerical results for each lattice geometry analyzed, Z–Axis mechanical study.

The methodology proposed and developed in this manuscript, as well as the set of numerical simulations performed have been carried out using a computer equipment of the MSI notebook type with an Intel (R) Core (TM) i-77700HQ CPU @ 2.80 GHz. On the other hand, regarding the computational cost, the time spent on the complete execution (preprocessing, processing and postprocessing phase) is estimated at 15 minutes.

## 4 Conclusions

In the present manuscript, a novel lattice-type structure with a pyramidal geometry is designed for 3D additive manufacturing processes. The proposed lattice unit cell demonstrates enhanced mechanical rigidity owing to its topology, which efficiently distributes the applied forces and resists effectively the obtained deformation. The parameterization process, based on Bounding Box dimensions, aspect ratios and lattice thickness, allows for versatile customization of the lattice unit cell to meet specific mechanical and geometric constraints.

An analytical evaluation of the structural behavior of the lattice unit cells validates the proposed design and parameterization process. Both the analytical and numerical approaches yield consistent results, confirming the structural integrity and applicability of the proposed lattice structure. The analytical model employed, rooted in the rigidity model of unidirectional 3D beam elements, provides a robust framework for understanding the structural behavior of the lattice unit cell.

In a case study involving a plastic part with defined bounding box dimensions and two commonly used lattice patterns (cross and hexagonal), the proposed lattice structure demonstrated superior structural efficiency.

The structural efficiency, defined as the ratio between the maximum displacements and mass, highlighted the optimized use of the material in the proposed lattice design, making it a promising candidate for lightweight and high-performance applications. Furthermore, the relative consistency in the structural efficiency across different load scenarios and lattice geometries emphasized the reliability and versatility of the proposed lattice structure. It offers an efficient balance between structural rigidity and material utilization, making it suitable for various engineering applications in which strength-to-weight ratio optimization is critical. Furthermore, in load states or mechanical stresses of greater relevance, where the magnitude of the displacement field is greater. In this way, it has been proven that the design of the proposed lattice unit cell, for the case study analyzed, improves the magnitude of the structural efficiency parameter for the load scenarios in the Y and Z direction, in which the magnitude of the displacement field is greater. In addition to the FEM analysis for structural performance, the manufacturability of the proposed lattice unit cell was simulated using additive manufacturing software. These simulations confirmed that the lattice design does not require support structures during printing, and it can be manufactured without difficulty using 3D additive manufacturing technology. The design was optimized to eliminate overhangs and complex geometries that typically necessitate support materials, ensuring a straightforward and efficient printing process.

In summary, the developed lattice structure, characterized by its pyramidal geometry and adaptive parameterization, presents a viable and efficient solution for enhancing the mechanical performance of components manufactured through 3D additive manufacturing. However, the manufacturability of such a lattice structure requires careful consideration of several factors inherent to 3D additive manufacturing processes. Challenges such as retraction, beam deflection due to thermal effects, and material deposition accuracy must be addressed, particularly for fine lattice geometries with small cross-sectional areas.

The design of the proposed lattice accounts for these challenges by incorporating relatively larger cross-sections and a simplified, support-free geometry to reduce manufacturing complexity. Future work may explore further optimizations to address these issues, ensuring that the lattice structure can be efficiently and reliably produced using 3D additive manufacturing technology, focusing on exploring additional lattice geometries, material variations, and applications in other areas of engineering and industrial design to further validate and extend the findings of the proposed design and developed methodology.

**Funding**

The research has been funded by the University of Jaén through the Research Plan 2022–2023-ACCION1a POAI 2022–2023: TIC-159, and by the scientific society INGEGRAF through the award to the authors for their research presented in the field of industrial design (2024).